\documentclass[a4paper,11pt]{article}
\usepackage{ifthen,calc}
\usepackage{amsmath,amssymb,amsthm,amstext,amsfonts,bm,mathrsfs,commath,comment}
\usepackage[scale=0.76]{geometry}
\usepackage{color,graphicx,pgf,tikz}
\usepackage{placeins}	% control float placement \FloatBarrier
\usepackage{float}
\usepackage{subfigure}
\usepackage{extarrows}
\usepackage{natbib}
\usepackage{enumitem}
\setlist{nosep}

\usepackage{setspace}
\allowdisplaybreaks
\usepackage{lineno}

\newtheorem{defn}{Definition} 
\newtheorem{thm}{Theorem}
\newtheorem{lem}{Lemma}

\newtheorem{prop}{Proposition}
\newtheorem{rmk}{Remark}
\newtheorem{exam}{Example}
\newtheorem{cro}{Corollary}

\newtheorem{claim}{Claim}

\newcommand{\bE}{\mathbb{E}}
\newcommand{\bN}{\mathbb{N}}
\newcommand{\bR}{\mathbb{R}}
\newcommand{\bP}{\mathbb{P}}
\newcommand{\bZ}{\mathbb{Z}}
\newcommand{\cB}{\mathcal{B}}
\newcommand{\cA}{\mathcal{A}}

\newcommand{\cC}{\mathcal{C}}

\newcommand{\cF}{\mathcal{F}}

\newcommand{\cI}{\mathcal{I}}

\newcommand{\cM}{\mathcal{M}}

\newcommand{\cS}{\mathcal{S}}

\newcommand{\cU}{\mathcal{U}}

\newcommand{\tens}[1]{%
  \mathbin{\mathop{\otimes}\limits_{#1}}%
}

\DeclareMathOperator{\rmd}{d\!}

\usepackage{hyperref}
\hypersetup{colorlinks=true,
	citecolor=blue,
	linkcolor=blue,
	urlcolor=blue,
	bookmarksopen=true,
	pdfstartview={XYZ null null 1.50},
	pdfpagelayout={SinglePage}
}

\usepackage{scalerel,stackengine}
\stackMath
\newcommand\ehat[1]{%
	\savestack{\tmpbox}{\stretchto{%
			\scaleto{%
				\scalerel*[\widthof{\ensuremath{#1}}]{\kern-.6pt\bigwedge\kern-.6pt}%
				{\rule[-\textheight/2]{1ex}{\textheight}}%WIDTH-LIMITED BIG WEDGE
			}{\textheight}%
		}{0.5ex}}%
	\stackon[1pt]{#1}{\tmpbox}%
}

\usetikzlibrary{calc,trees,positioning,arrows,chains,shapes,matrix,automata,patterns,shapes.geometric,decorations.pathreplacing,decorations.pathmorphing,decorations.markings,shapes.symbols}
\tikzset{>=latex}

\allowdisplaybreaks[4]

\begin{document}
\title{Equilibrium convergence in large games\thanks{The first version of this paper was titled ``A unifying closed graph property for equilibria in large games'' in November 2020 (\href{https://arxiv.org/abs/2011.06789v1}{https://arxiv.org/abs/2011.06789v1}). Part of the results in this paper were presented at the 21st SAET Conference, ANU, Canberra, Australia, July 16--22, 2022, the NUS microeconomic theory seminar, June 17, 2022; and we thank the participants for their constructive comments and suggestions. Enxian~Chen's research is supported by the National Natural Science Foundation of China (No.~72303115) and the Fundamental Research Funds for the Central Universities (No.~63222003). Bin Wu acknowledges the National Natural Science Foundation of China (No.~72203154, No.~72450003).}}

\author{Enxian~Chen\thanks{School of Economics, Nankai University, Tianjin, 300071, China. E-mail: \href{mailto:chenenxian777@gmail.com}{chenenxian777@gmail.com}.}
\and
Bin~Wu\thanks{International School of Economics and Management, Capital University of Economics and Business, Fengtai District, Beijing, 100070, China. E-mail: \href{bwu.econ@gmail.com}{bwu.econ@gmail.com.}}
\and
Hanping Xu\thanks{Department of Economics, Chinese University of Hong Kong, Shatin, NT, Hong Kong SAR, 999077, China. E-mail: \href{hanpingxu@cuhk.edu.hk}{hanpingxu@cuhk.edu.hk}.}
}

\date{This version: \today}

\maketitle
\thispagestyle{empty} 

\begin{abstract}
This paper presents a general closed graph property for (randomized strategy) Nash equilibrium correspondence in large games. In particular, we show that for any large game with a convergent sequence of finite-player games, the limit of any convergent sequence of Nash equilibria of the corresponding finite-player games can be induced by a Nash equilibrium of the large game. 
Such a result goes beyond earlier results on the closed graph property for pure strategy Nash equilibrium correspondence in large games in multiple aspects. An application on equilibrium selection in large games is also presented.

\bigskip
\textbf{JEL classification}: C60; C62; C72

\bigskip
\textbf{Keywords}:  Large games, finite-player games, Nash equilibrium, closed graph property.

\end{abstract}

\clearpage
\tableofcontents
%\setlength{\parskip}{4pt}
%\linenumbers	% test command  
\clearpage

\newpage  
\setcounter{page}{1}
\section{Introduction}\label{sec-intro}
Games with a continuum of agents (henceforth referred to as ``large games'') have been extensively studied in the literature. These games model strategic interactions involving a large number of individuals, each with a negligible ability to affect the others. Theoretically, large games exhibit various desirable properties, such as the existence of a pure strategy Nash equilibrium. In practical applications, the model of large games is also widely used to describe real-life situations involving many agents, such as in matching, contests, finance, and so on.

Given that large games involve infinitely many players while real-life situations involve only finitely many, it is important to examine the relationship between the Nash equilibria of a large game and those of a convergent sequence of games with large but finitely many players. Specifically, a fundamental question in the theory of large games is whether any sequence of Nash equilibria of a sequence of finite-player games converging to a limit game, converges to a Nash equilibrium of the limit game. This property, known as the closed graph property of the Nash equilibrium correspondence, has been extensively studied in the context of pure strategies.\footnote{See, for example,  \cite{KRSY2013}, \cite{QY2014}, \cite{QYZ2016}, and \cite{HSS2017}.} However, focusing solely on pure strategy Nash equilibria reveals various limitations in considering the closed graph property.

Firstly, it is well-known that finite-player games may not always have a pure strategy Nash equilibrium. For instance, Example \ref{exam-RSP} in Subsection~\ref{subsec-motivation} presents a convergent sequence of finite-player games that lack pure strategy Nash equilibria, generalizing the classical two-player Rock-Scissors-Paper  game to many players. Thus, earlier results on the closed graph property of Nash equilibrium correspondence are not applicable in this case due to the absence of pure strategy Nash equilibria. 

Secondly, even if the closed graph property holds for a convergent sequence of finite-player games with convergent pure strategy Nash equilibria, it does not address convergent sequences of randomized (non-pure) strategy Nash equilibria for the underlying sequence of finite-player games. For instance, Example \ref{exam-invest} in Subsection~\ref{subsec-motivation} illustrates a convergent sequence of investment games that have pure strategy Nash equilibria. However, the limit of this sequence leads to a large game with a natural Nash equilibrium that is not a limit point of any sequence of pure strategy Nash equilibria from the finite-player games; rather, it is a limit point of the randomized strategy Nash equilibria of those finite-player games. Thus, focusing solely on pure strategy Nash equilibria limits the applicability of the closed graph property.

Thirdly, randomized strategy Nash equilibrium are widely used and studied in many applied large economy models. For example, \cite{S2006} considered a large market where a group of $n$ ``quacks'' engage in price competition over a continuum of patients. Assuming that patients behave according to a boundedly rational procedure, \cite{S2006} demonstrated that for any $n \ge 2$, the game has a unique Nash equilibrium in which every quack employs a (symmetric) randomized strategy.

Therefore, to address these issues in both theoretical and applied models, we need to consider a more general version of the closed graph property that includes randomized strategy profiles. As the main result of this paper, we show that for any large game with a convergent sequence of finite-player games, the limit distribution of any convergent sequence of (randomized strategy) Nash equilibria from the corresponding finite-player games can be induced by a Nash equilibrium in the limit large game. This general closed graph property addresses all the aforementioned limitations.

In particular, our main result extends beyond earlier findings on the closed graph property for pure strategy Nash equilibria in the following ways: (i) Since randomized strategy Nash equilibria always exist in both finite-player and large games, our result applies to any large game and its finite approximations.\footnote{The multi-player Rock-Scissors-Paper game and the investment game discussed in Subsection~\ref{subsec-motivation} feature randomized strategy Nash equilibria where all players select each action with equal probability. It is clear that this sequence of Nash equilibria converges to a randomized strategy Nash equilibrium in the limit large game.} (ii) Our result encompasses any convergent sequence of Nash equilibria in finite-player games and asserts that the limit is a Nash equilibrium. This provides a comprehensive answer to the closed graph property of Nash equilibrium correspondence in large games. In addition, unlike earlier proofs in the literature for the case of pure strategy equilibria, a technical difficulty arises in our proof of the main result. Namely, to show the convergence of Nash equilibria, we need to estimate the gap between the (expected) equilibrium payoff and the realized payoffs of a (randomized strategy) equilibrium.\footnote{Such an issue does not occur when dealing with the convergence of pure strategy Nash equilibria. While some papers in the literature also examine randomized strategy profiles/Nash equilibria, they address different research questions. For more details, see Section~\ref{sec-literature}.}

Next, we discuss two extensions of the main theorem. By allowing the games in the convergent sequence to be large games, we introduce the concept of the strong closed graph property and demonstrate that the Nash equilibrium correspondence of any large game satisfies this property. Furthermore, we show that under certain additional conditions on the large game, the Nash equilibrium correspondence also exhibits the pure closed graph property. Specifically, we prove that for any large game with a convergent sequence of finite-player games, the limit of any convergent sequence of Nash equilibria from the corresponding finite-player games can be induced by a pure strategy Nash equilibrium of the large game.

As an application of our main result, we introduce an equilibrium selection method for large games. Equilibrium selection is an important topic in game theory because games typically have multiple Nash equilibria, some of which may be unreliable and need to be excluded. Various refined solution concepts have been proposed in the literature based on different selection criteria. In Section~\ref{sec-application}, we present a new equilibrium selection criterion for large games. The main idea is as follows: Suppose a group of $n$ players is involved in a large finite-player game $G^n$. To simplify the analysis, they use a large game $G$ with a continuum of players for approximation.  In many cases, $G$ has multiple Nash equilibria, denoted by $E$. We can divide $E$ into two disjoint subsets:  $E_1$ contains equilibria of $G$ that are limit points of a sequence of equilibria of $\{G^n\}_{n\in \bZ_+}$, while $E_2$ consists of non-approximable equilibria. Since players are in fact playing the large finite game $G^n$, only the equilibria in $E_1$ are considered reliable, while those equilibria in $E_2$ are ``singular'' points generated during the idealization process and should be excluded. In Corollary~\ref{application}, we show that $E_1$ is non-empty,  ensuring that an approximable equilibrium always exists in a large game with a sequence of finite approximations.

The remainder of this paper is structured as follows: Section~\ref{sec-model} introduces the models of large games and finite-player games, along with the concept of Nash equilibrium in these contexts. Our main results are presented in Section~\ref{sec-main results}. Section~\ref{sec-extension} discusses two extensions of the main theorem. An application is outlined in Section~\ref{sec-application}. Section~\ref{sec-literature} reviews related literature, while the proofs of our results are provided in Section~\ref{sec-appendix}.

%-----------------------------------------------------------------------------------

\section{Basic models}\label{sec-model}
In this section, we introduce basic definitions of large games and finite-player games. In these games, each player has a general compact set of actions, and her payoff depends on her own choice as well as  the distribution of actions chosen by all players. Specifically,  our model allows different players to have different feasible action sets.

\subsection{Large games}\label{subsec-large}
A large game is defined as follows. Let $(I, \cI, \lambda)$ be an atomless probability space representing the set of players,\footnote{Throughout this paper, we follow the convention that a probability space is complete and countably additive.} and let $A$ be a compact metric space representing a common action space, endowed with the Borel $\sigma$-algebra $\cB(A)$. Let $\cA \colon I \twoheadrightarrow A$ be a nonempty, measurable, and compact-valued correspondence, specifying a feasible action set  $A_i=\cA(i)$ for each player $i\in I$. The set of Borel probability measures on $A$ is denoted by $\cM(A)$. Given an action profile of all the players, the action distribution, which specifies the portions of players taking some actions in $A$ (also called a societal summary), can be viewed as an element in $\cM(A)$. Each player's payoff is a bounded continuous function on $A \times \cM(A)$, which means that the payoff continuously depends on player's own choice and the societal summaries. Let $\cU_A$ be the space of bounded continuous functions on $A \times \cM(A)$, endowed with the sup-norm topology and the resulting Borel $\sigma$-algebra.

Let $\cC_A$ be the set of all compact subsets of $A$. Endowed with the Hausdorff metric and the Borel $\sigma$-algebra, $\cC_A$ is a compact metric space (\citet[Theorem 3.85]{AB2006}). For each player $i \in I$, their characteristic consists a feasible action set $A_i$ (an element of $\cC_A$) and a payoff function $u_i$ (an element of $\cU_A$). Thus, the space of all players' characteristics is $\cC_A \times \mathcal{U}_A$, endowed with the product topology.

A \emph{large game} $G$ is a measurable mapping from $(I, \cI, \lambda)$ to $\cC_A \times \mathcal{U}_A$. A \emph{pure strategy profile} $f$ is a measurable mapping from $(I, \cI, \lambda)$ to $A$ such that for $\lambda$-almost all $i \in I$, $f(i) \in A_i$. Let $\lambda f^{-1}$ be the societal summary (also denoted by $s(f)$), which represents the societal action distribution induced by $f$. Specifically, $s(f) (B)$,  where $B$ is a subset of $A$, denotes the portion of players taking actions in $B$.

A randomized strategy for player $i$ is a probability distribution $\mu \in \cM(A_i)$. A \emph{randomized strategy profile} $g$ is a measurable function from $I$ to $\cM(A)$ such that $g(i, A_i) = 1$ for $\lambda$-almost all $i \in I$. Notice that every pure strategy profile $f$ naturally corresponds to a randomized strategy profile $g^f$ where $g^f(i) = \delta_{f(i)}$\footnote{Here $\delta_{f(i)}$ denotes the Dirac probability measure that assigns probability one to $\{f(i)\}$.} for each player $i \in I$. Given a randomized strategy profile $g$, we model the societal summary $s(g)$ as the average action distribution of all the players, i.e., $ s(g) = \int_I g(i) \dif\lambda(i) \in \cM(A)$.\footnote{Note that the societal summary $\int_I  g(i) {\rm{d}}\lambda(i)$ is an element in $\cM(A)$ that satisfies $\int_I  g(i) {\rm{d}}\lambda(i)(B)=\int_I  g(i, B) {\rm{d}}\lambda(i)$
for all $B\in\cB(A).$} Clearly, when $f$ is a pure strategy profile, $\int_I f(i) \dif\lambda(i)$ reduces to $\lambda f^{-1}$, which is the societal action distribution induced by $f$.

The formal definition of a randomized strategy Nash equilibrium is stated as follows.
\begin{defn}[Randomized strategy Nash equilibrium]
\label{defn:randNE}
\rm
A randomized strategy profile $g \colon I \rightarrow \cM(A)$ is said to be a \emph{randomized strategy Nash equilibrium} if for $\lambda$-almost all $i \in I$,
$$ \int_{A_i} u_i \bigl( a, s(g) \bigr) g(i, \dif a)  \ge  \int_{A_i} u_i \bigl( a, s(g) \bigr) \dif\mu(a) \text{ for all } \mu \in \cM(A_i). $$
\end{defn}

A randomized strategy profile $g$ is a randomized strategy Nash equilibrium if it is optimal for almost all players with respect to the societal summary $s(g)$ in terms of expected payoff. For a pure strategy profile $f$, its average action distribution $s(f)$ reduces to the societal action distribution $\lambda f^{-1}$. Thus, we have the following definition of pure strategy Nash equilibrium.

\begin{defn}[Pure strategy Nash equilibrium]
\label{defn:pureNE}
\rm
A pure strategy profile $f \colon I \rightarrow A$ is said to be a \emph{pure strategy Nash equilibrium} if, for $\lambda$-almost all $i\in I,$
$$u_{i}\bigl(f(i),  \lambda f^{-1}\bigr) \ge u_{i}(a, \lambda f^{-1}) \text{ for all } a \in A_i.$$
\end{defn}

\begin{rmk}
\rm
Our model of large games is general enough to encompass the setting of large games with traits studied in the literature (see, for example, \cite{KRSY2013}).\footnote{Let $T$ be a compact metric space representing the space of traits, endowed with its Borel $\sigma$-algebra. The trait function $\alpha \colon I \to T$ is defined as a measurable mapping that associates each player with a trait. Let $\cU_{(A,T)}$ be the space of bounded and continuous real-valued functions on $A \times \cM(T \times A)$, where $\cM(T \times A)$ is the set of Borel probability measures on $T \times A$.  A large game with traits $G$ is a measurable mapping from $I$ to $T \times \cU_{(A,T)}$ such that $G = ( \alpha, u)$. Then we consider another large game without traits $\widetilde{G}$ as follows. Let the action set $\widetilde{A} = T \times A$, and the action correspondence $\widetilde{\cA}(i) = \{\alpha(i) \} \times A_i$. Let player $i$'s payoff function $\widetilde{u}_i$ be defined as $\widetilde{u}_i(t, a, \tau) = u_i(a, \tau)$ for any $t\in T$, $a \in A_i$, and $\tau \in \cM(T \times A)$. Since $(t, a, \tau) \in  T \times A\times \cM( T \times A) = \widetilde{A} \times \cM(\widetilde{A})$, $\widetilde{u}_i$ is a bounded continuous function on $\widetilde{A} \times \cM(\widetilde{A})$. Let the large game $\widetilde{G}\colon I \to \cC_{\widetilde{A}} \times \cU_{\widetilde{A}}$ be defined as $\widetilde{G}(i) = (\widetilde{\cA}(i), \widetilde{u}_i)$. Thus, a large game $G$ with traits can be viewed as a large game $\widetilde{G}$ without traits.} %We implicitly assume that players make independent choices within this model. However, this independence assumption is incompatible with the measurability of the strategy profile. To address this,  we adopt the rich Fubini extension framework  introduced in \cite{Sun2006} and assume that the player space $(I, \cI, \lambda)$ allows such an extension.\footnote{See \citet[Section 5]{Sun2006}, \cite{SZ2009}, and \cite{Podczeck2010} for constructions of rich Fubini extensions.} Moreover, a mixed strategy Nash equilibrium on a Fubini extension space is equivalent to a randomized strategy Nash equilibrium as defined in Definition~\ref{defn:randNE} (\cite{KRSY2015}).
\end{rmk}

\subsection{Finite-player games}\label{subsec-finite}

In this subsection, we introduce a class of finite-player games where each player's payoff function depends on her own choice and the societal summary. Let  $(I^n, \cI^n, \lambda^n)$ denote  the finite probability space representing the set of players. Here, $|I^n| = n$ and $\cI^n$ consists of all the subsets of $I^n$ (i.e., the power set of $I^n$). For each player $i \in I^n$, her action set is $A_i^n$, a nonempty and closed subset of the common compact action space $A$. Similarly, the correspondence $\cA^n \colon I^n \twoheadrightarrow A$ that satisfies $\cA^n(i) = A^n_i$ as a representation of the action correspondence. Each player's payoff function depends on her own choice and the action distribution induced by the choices of all players in $I^n$ (i.e., the societal summary). Clearly, the set of such action distributions is a subset of $\cM(A)$ and is denoted by
$$
D^n = \Bigl\{   \tau \in \cM(A) \;\Big\vert\; \tau =  \sum\limits_{i \in I^n} \lambda^n(i) \delta_{ a_i }  \text{  where } a_i \in A^n_i \text{ for all } i \in I^n          \Bigr\}.
$$
Player $i$'s payoff function is then given by a bounded continuous function $u_i^n \colon A \times \cM(A) \to \bR$. Thus, a finite-player game $G^n$ can be viewed as a mapping from $I^n$ to $\cC_A \times \mathcal{U}_A$ such that $G(i) = (A_i^n, u_i^n)$ for all $i \in I^n$.

In this finite-player game, a \emph{pure strategy profile} $f^n$ is a mapping from $( I^n , \cI^n, \lambda^n)$ to $A$ such that $f^n( i ) \in A^n_i$ for all $i \in I^n$. Hence, given a pure strategy profile $f^n$, the payoff function for player $i$ is
$$
u_i^n(f^n) = u_i^n \Bigl( f^n(i), \sum\limits_{j \in I^n } \lambda^n(j) \delta_{f^n(j)} \Bigr),
$$
here we slightly abuse the notation $u_i^n(f^n)$ to denote player $i$'s payoff given the strategy profile $f^n$.

Similarly, a randomized strategy of player $i$ is a probability distribution $\mu \in \cM(A^n_i)$. A \emph{randomized strategy profile} $g^n$ is a mapping from  $( I^n , \cI^n, \lambda^n)$ to $\cM(A)$ such that $g^n(i , A^n_i) = 1$. Thus, given a randomized strategy profile $g^n$, player $i$'s (expected) payoff is
$$
u_i^n(g^n) = \int_{\prod\limits_{j \in I^n} A^n_j} u_i^n \Bigl(a_i, \sum\limits_{j \in I^n } \lambda^n(j) \delta_{a_j} \Bigr) \tens{j \in I^n} g^n(j, {\rm{d}} a_j),
$$
where $\tens{ j \in I^n} g^n(j,  {\rm{d}} a_j)$ is the product probability measure on the product space $\prod\limits_{j \in I^n} A^n_j$. The societal summary induced by $g^n$ is $s(g^n) = \int_{I^n} g^n(i) \mathrm{d} \lambda^n(i)$. Finally, we define a randomized strategy Nash equilibrium as follows.

\begin{defn}[Randomized strategy Nash equilibrium]
	\label{defn:randNE-finite}
	\rm
	A randomized strategy profile $g^n \colon I^n \rightarrow \cM(A)$ is said to be a \emph{randomized strategy Nash equilibrium} if for all $i \in I^n$,
	$$ u_i^n(g^n) \ge  u_i^n( \mu  ,g^n_{- i}  ) \text{ for all } \mu \in \cM(A^n_i), $$
	where $(\mu, g^n_{-i})$ represents the randomized strategy profile such that player $i$ plays the randomized strategy $\mu$, and player $j$ plays the randomized strategy $g^n(j)$ for all $j \in I^n \backslash \{i\}$.
\end{defn}

Throughout the rest of this paper, a Nash equilibrium always refers to a randomized strategy Nash equilibrium.

%-------------------------------------------------------------------------------------
\section{Main results}\label{sec-main results}
We are now prepared to present our main result. In Subsection~\ref{subsec-motivation}, we provide two motivating examples that illustrate the importance of considering randomized strategy profiles when examining equilibrium convergence in large games. Next, we introduce the concept of the general closed graph property and establish our main theorem in Subsection~\ref{subsec-closed}, demonstrating that the Nash equilibrium correspondence satisfies this property.

\subsection{Motivating examples}\label{subsec-motivation}
We first present an example that generalizes the classical two-player Rock-Paper-Scissors game to an $n$-player setting. As we will demonstrate, this game does not have any pure strategy Nash equilibrium when $n \ge 2$. As the number of players approaches infinity, the sequence of finite-player games converges to a large game played by a continuum of players. This large game features a unique equilibrium action distribution, where each action is chosen by one-third of the players. Additionally, this equilibrium action distribution can be approximated by a sequence of randomized strategy Nash equilibria from the finite-player games.

\begin{exam}\label{exam-RSP}
Fix an integer $n\geq 2$, we consider a Rock-Scissors-Paper  game $G^n$ with $n$ players as follows. Let $I^n$ denote the set of players and $|I^n| = n$ and let $\lambda_n$ be the counting measure on $I_n$. Each player must choose one of three actions: Rock ($R$), Scissors ($S$), or Paper ($P$), creating a common action set $A = \{  R,S,P \}$. Similar to the two-player game setting, a player who chooses Rock will beat another player who chooses Scissors but will lose to one who chooses Paper; a player who chooses Paper will lose to a player who chooses Scissors. For each player $i \in I^n$, her payoff is the difference between the portion of players she beats and the portion of players she loses to.\footnote{For example, suppose $n=3$ and players choose $R$, $S$, $P$ respectively. Then each player beats one of the rest two players, but loses to the other player, hence each player's payoff is $\tfrac{1}{3} - \tfrac{1}{3} = 0$.} Thus, each player's payoff function is given as follows:
\begin{align*}
		&  u(a, \tau) =
		\left\{
		\begin{array}{rcl}
			\tau(S) - \tau(P)    &   & {\textrm{if}  \ a = R,}\\
			\tau(P) - \tau(R)    &   & {\textrm{if}  \ a = S,}\\
			\tau(R) - \tau(S)    &   & {\textrm{if}  \ a = P,}\\
		\end{array}
		\right.
	\end{align*}
where $\tau(a)$ denotes the proportion of players choosing the action $a \in \{R, S, P\}$.
\end{exam}

Notice that the game is zero-sum and when $n=2$, it is exactly the classical two-player Rock-Scissors-Paper  game.

\begin{claim}\label{claim-no PSNE}
There is no pure strategy Nash equilibrium in the game $G^n$ when $n \ge 2$.
\end{claim}

Clearly, if the number of players tends to infinity, the sequence of finite-player games $\{G^n\}_{n \in \bZ_+}$ converges to a large Rock-Scissors-Paper  game played by a continuum of players. 
Hence, Claim~\ref{claim-no PSNE} highlights the importance of considering randomized strategy profiles when using a sequence of finite-player games to approximate the large game.

Let $G$ denote the limit large game, and let $(I,\mathcal{I}, \lambda)$ denote the player space, where $I$ is the Lebesgue unit interval . All the players have the same action set $A$ and the same payoff function $u$ as in finite-player game $G^n$. Since the player space $I$ is atomless, each player in the large game $G$ has no influence on the aggregate action distribution. We will  see that $G$ has a unique Nash equilibrium action distribution where each action is chosen by one-third of the players.

\begin{claim}\label{claim-large symmetric}
The large Rock-Scissors-Paper  game $G$ has a unique Nash equilibrium action distribution $\tau^* = (\tfrac{1}{3}, \tfrac{1}{3}, \tfrac{1}{3})$.
\end{claim}

Note that the equilibrium action distribution $\tau^* = (\tfrac{1}{3}, \tfrac{1}{3}, \tfrac{1}{3})$ can be implemented by a symmetric Nash equilibrium of the large game $G$. That is, players can simply coordinate on the same strategy $\frac{1}{3}\delta_{R} + \frac{1}{3}\delta_{S} + \frac{1}{3}\delta_{P}$ (i.e., choose each action with equal probability). As shown by the following result, this strategy profile is also a Nash equilibrium in each game $G^n$, and hence the unique equilibrium action distribution of the large game $G$ can be approximated by a sequence of randomized strategy Nash equilibria of finite-player games $\{G^n\}_{n \in \bZ_+}$.

\begin{claim}\label{claim-finite symmetric}
For each $n \ge 2$, the  strategy profile $g^n(i) \equiv \frac{1}{3}\delta_{R} + \frac{1}{3}\delta_{S} + \frac{1}{3}\delta_{P}$ is a Nash equilibrium of the game $G^n$.
\end{claim}

\bigskip

The sequence of finite-player games $\{G^n\}_{n \in \bZ_+}$ discussed in the previous example has no pure strategy Nash equilibrium. Below we present another example of finite-player games $\{G^n\}_{n \in \bZ_+}$ that do have pure strategy Nash equilibria. However, the limit large game features a Nash equilibrium (action distribution) that is a limit point of randomized (non-pure) strategy Nash equilibria of $\{G^n\}_{n \in \bZ_+}$.

\begin{exam}\label{exam-invest}
A group of $n$ traders in a financial market needs to make investment decisions involving two similar shares, $a_1$ and $a_2$. Let $I^n$ denote the set of traders and $|I^n| = n$ and let $\lambda_n$ be the counting measure on $I_n$. Each trader's payoff depends on her own choice as well as the proportion of traders making the same choice, given by
\begin{align*}
		&  u(a, \tau) =
		\left\{
		\begin{array}{rcl}
			h\bigl(\tau(a_1)\bigr)    &   & {\textrm{if}  \ a = a_1,}\\
			h\bigl(\tau(a_2)\bigr)    &   & {\textrm{if}  \ a = a_2,}\\
		\end{array}
		\right.
	\end{align*}
where $h \colon [0, 1] \to \bR$ is a strictly increasing and continuous function. Let $G^n$ denote this game and $G$ represent the limit large game with a continuum of players.
\end{exam}

Clearly, $G^n$ have pure strategy Nash equilibria for each $n \ge 2$; for example, all players can choose $a_1$ (or $a_2$) simultaneously. On the other hand, the limit large game $G$ has a natural Nash equilibrium where $a_1$ and $a_2$ are chosen by exactly half of the players. The induced equilibrium action distribution is $\tau^* = (\tfrac{1}{2}, \tfrac{1}{2})$.  

\begin{claim}\label{claim-invest}
$\tau^*$ is not a limit point of any sequence of pure strategy Nash equilibria of $\{G^n\}_{n \in \bZ_+}$.  
\end{claim}
Notice that the convergence refers to the convergence of equilibrium action distributions. Claim~\ref{claim-invest} does not imply that $G$ is not an idealized approximation of $G^n$. In fact, it is straightforward to verify that the randomized strategy profile $g^n(i) \equiv \frac{1}{2}\delta_{a_1} + \frac{1}{2}\delta_{a_2}$ is a Nash equilibrium of $G^n$ for each $n \ge 2$. Furthermore, it is clear that $\tau^*$ is the limit point of $\{g^n\}_{n \in \bZ_+}$.

\subsection{General closed graph property}\label{subsec-closed}
In this subsection, we establish a theory regarding the closed graph property for the Nash equilibrium correspondence of the large game. We will show that for any sequence of large finite games  converging to the limit game  and any sequence of Nash equilibria corresponding to these finite games, the weak limit of the sequence of Nash equilibria must be induced by a Nash equilibrium of the limit game.

Let $\{( I^{n},\cI^{n}, \lambda^{n})\}_{n\in \bZ_{+}}$ be a sequence of probability spaces where $|I^n| = n$, $\mathcal{I}^{n}$ is the power set of $ I^{n}$, and $\lambda^n$ is a probability measure on $I^n$ such that $\sup_{i\in I^{n}} \lambda^{n}( i)\rightarrow 0$ as $n$ goes to infinity. For each $n\in {\mathbb{Z}}_{+},$ let a finite-player game $G^{n}=(\cA^{n}, u^{n})$ be a mapping from the player space $(I^{n},\mathcal{I}^{n}, \lambda^{n})$ to the characteristic space $\cC_A \times \mathcal{U}_A$. We state the formal definition of the general closed graph property as follows.

\begin{defn}[General closed graph property]
\label{defn-closed graph property}
\rm
The  Nash equilibrium correspondence of a large game $G \colon I \to \cC_A \times \mathcal{U}_A$ is said to have the \emph{general closed graph property} if
\begin{itemize}
	\item[(i)] for any sequence of finite-player games $\{G^n\}_{n \in \bZ_{+}}$ converging to the large game $G$ in the sense that $\{\lambda^n (G^n)^{-1} \}_{n \in \bZ_{+}}$ converges weakly to $\lambda G^{-1}$, and
	\item[(ii)] for any sequence $\{g^n\}_{n \in \bZ_{+}}$ where each $g^n$ is a Nash equilibrium of $G^n$ such that the sequence of societal summaries $\{s(g^n)\}_{n \in \bZ_{+}}$ converges weakly to a distribution $\tau^*$ on $A$,
\end{itemize}
then there exists a Nash equilibrium $g$ of $G$ such that $s(g) = \tau^*$.
\end{defn}

Notice that the societal summary $s(g^n) = \int_{I^n} g^n(i) \mathrm{d} \lambda^n(i)$ is induced by the randomized strategy profile $g^n$ of the finite-player game $G^n$. Since $G^n$ is a measurable mapping from $I^n$ to $\cC_A \times \mathcal{U}_A$, $\lambda^n (G^n)^{-1}$ is an induced measure on $\cC_A \times \mathcal{U}_A$. Similarly, $\lambda G^{-1}$ is also an induced measure on $\cC_A \times \mathcal{U}_A$. Thus, the convergence of games is defined as the weak convergence of measures on $\cC_A \times \mathcal{U}_A$.

This definition of the closed graph property generalizes the definition in the literature that focuses on the pure strategy Nash equilibrium (See, for example, \cite{KRSY2013}, \cite{QY2014}, \cite{QYZ2016}, and \cite{HSS2017}). Based on our discussions of the motivating examples in the previous subsection, it is essential to extend the closed graph property to include randomized strategy Nash equilibrium correspondences. We are now ready to present the main result of this paper as follows.

\begin{thm}\label{thm-closed graph property}
	The Nash equilibrium correspondence of any large game $G$ has the general closed graph property.
\end{thm}

Theorem~\ref{thm-closed graph property} extends existing results on the closed graph property by demonstrating that any convergent sequence of Nash equilibria from finite-player games converges to a Nash equilibrium of the limit large game. The proof of Theorem~\ref{thm-closed graph property} is more complex than those proofs of the closed graph property in the case of pure strategies. The primary difficulty lies in considering all  possible realizations of mixed Nash equilibria of finite-player games and their convergence. The detailed proof is in Subsection~\ref{subsec-proof main results}.

%-------------------------------------------------------------------------

\section{Extensions}\label{sec-extension}
In this section, we explore two extensions of the main theorem. Subsection~\ref{subsec-strong} introduces the concept of the strong closed graph property, which accommodates the presence of large games in the convergent sequence. We then examine another refined version of the general closed graph property in Subsection~\ref{subsec-pure}, termed the pure closed graph property, which requires the limit Nash equilibrium to be in pure strategies.

\subsection{Strong closed graph property}\label{subsec-strong}
Definition~\ref{defn-closed graph property} of the general closed graph property can be extended to accommodate  for ``non-atomic'' games in the convergent sequence. Specifically, the convergent sequence of games $\{G^n\}_{n \in \bZ_{+}}$ may include large games with a continuum of players. This stronger property is referred to as the strong closed graph property, defined as follows.

\begin{defn}[Strong closed graph property]
\label{defn-strong closed graph property}
\rm
The  Nash equilibrium correspondence of a large game $G \colon I \to \cC_A \times \mathcal{U}_A$ is said to have the \emph{strong closed graph property} if
\begin{itemize}
	\item[(i)] for any sequence of games $\{G^n\}_{n \in \bZ_{+}}$ converging to the large game $G$ in the sense that $\{\lambda^n (G^n)^{-1} \}_{n \in \bZ_{+}}$ converges weakly to $\lambda G^{-1}$, and
	\item[(ii)] for any sequence $\{g^n\}_{n \in \bZ_{+}}$ where each $g^n$ is a Nash equilibrium of $G^n$ such that the sequence of societal summaries $\{s(g^n)\}_{n \in \bZ_{+}}$ converges weakly to a distribution $\tau^*$ on $A$,
\end{itemize}
then there exists a Nash equilibrium $g$ of $G$ such that $s(g) = \tau^*$.
\end{defn}

\begin{prop}\label{prop-strong}
	The Nash equilibrium correspondence of any large game $G$ has the strong closed graph property.
\end{prop}

The detailed proof of Proposition~\ref{prop-strong} is provided in Subsection~\ref{subsec-proof-strong}. It is evident that we only need to consider the case where $\{G^n \}_{n \in \bZ_{+}}$ is a sequence of large games.

\subsection{Pure closed graph property}\label{subsec-pure}
The general closed graph property does not require that the limit Nash equilibrium of the large game to be a pure strategy profile. In this subsection, we consider another refinement of Definition~\ref{defn-closed graph property}, termed the pure closed graph property, defined as follows.

\begin{defn}[Pure closed graph property]
	\label{defn-strong closed graph property}
	\rm
	The  Nash equilibrium correspondence of a large game $G \colon I \to \cC_A \times \mathcal{U}_A$ is said to have  \emph{pure closed graph property} if
	\begin{itemize}
		\item[(i)] for any sequence of finite-player games $\{G^n\}_{n \in \bZ_{+}}$ converging to the large game $G$ in the sense that $\{\lambda^n (G^n)^{-1} \}_{n \in \bZ_{+}}$ converges weakly to $\lambda G^{-1}$, and
		\item[(ii)] for any sequence $\{g^n\}_{n \in \bZ_{+}}$ where each $g^n$ is a Nash equilibrium of $G^n$ such that the sequence of societal summaries $\{s(g^n)\}_{n \in \bZ_{+}}$ converges weakly to a distribution $\tau$ on $A$, 
	\end{itemize}
	then there exists a pure strategy Nash equilibrium $f$ of $G$ such that  $s(f) = \tau$.
\end{defn}

The pure closed graph property may not hold for some Nash equilibrium correspondences of large games. The reason is that a general large game my have no pure strategy Nash equilibrium; see, for instance, Example~1 of \cite{QY2014}. Therefore, to reestablish the pure closed graph property in large games, we need to impose additional conditions on the large game.

Existing literature has indicated that the existence of pure strategy Nash equilibria in large games may depend on the cardinality of the underlying set of actions. Specifically, if there are at most countably many actions, a pure strategy Nash equilibrium always exists. Additionally, another body of work has shown that pure strategy Nash equilibria exist in large games where the player space satisfies certain saturation conditions. Motivated by these findings in the literature, we impose the following two conditions.

\begin{defn}
	\label{defn-countable}
	\rm
	A large game $G$ is said to be \emph{countable-valued} if the action space $A$ is at most countable and the range of action correspondence $\cA \colon I \to \cC_A$ is at most countable.
\end{defn}

\begin{defn}
	\label{defn-nowhere}
	\rm
     Let $( I,\mathcal{I}, \lambda)$ be an atomless probability space, and $\cF$ be a sub-$\sigma$-algebra of $\cI$. The $\sigma$-algebra $\mathcal{I}$ is said to be \textit{nowhere equivalent} to its sub-$\sigma$-algebra $\mathcal{F}$ if for every nonnegligible subset $D \in \mathcal{I}$, there exists an $\mathcal{I}$-measurable subset $D_{0}$ of $D$ such that $ \lambda(D_{0}\triangle D_{1})>0$ for any $D_{1}\in \mathcal{F}^{D}$, where $D_{0}\triangle D_{1}$ is the symmetric difference $(D_0\setminus D_1)\cup (D_1\setminus D_0)$, and $\mathcal{F}^{D}$ is the restricted $\sigma$-algebra $\{D\cap D^{'}, D^{'}\in{\mathcal{F}}\}$.
\end{defn}

 For any large game $G \colon (I,\mathcal{I}, \lambda) \rightarrow \cC_A \times \mathcal{U}_A$, let $\sigma(G)$ be the  $\sigma$-algebra generated by $G$. That is, $\sigma(G)$ is the minimal $\sigma$-algebra of $I$ that makes $G$  measurable. Now we are ready to present the main result of this subsection. 

\begin{prop}\label{prop-refine}
A large game $G \colon ( I,\mathcal{I}, \lambda)\rightarrow \cC_A \times \mathcal{U}_A$, has the pure closed graph property if one of the following conditions holds:
\begin{itemize}
\item[(1)] $G$ is countable-valued;
\item[(2)]  $\mathcal{I}$ is nowhere equivalent to $\sigma(G)$.
\end{itemize}
\end{prop}

 \cite{QYZ2016} and \cite{HSS2017} studied the closed graph property in terms of pure strategies in large games under the condition of countable-valued large games and the condition of nowhere equivalence, respectively. Proposition~\ref{prop-refine} extends the results in \cite{QYZ2016} and \cite{HSS2017} in the following aspects: (i) it considers randomized strategy Nash equilibria within the convergent sequence of finite-player games; (ii) it accommodates the action correspondence in the large game model.

The proof of part (1) is a straightforward application of Theorem~\ref{thm-closed graph property} combined with the purification result from \citet[Theorem 2]{KRYZ2017}, so we omit the detailed proof in this paper. The detailed proof of part (2) is provided in Subsection~\ref{subsec-proof-pure}.\footnote{It is important to note that part (2) cannot be established by simply using Theorem~\ref{thm-closed graph property} to demonstrate that the limit distribution $\tau$ can be induced by a Nash equilibrium $g$, and then purifying a pure strategy profile $f$ from $g$ to conclude that $f$ is the desired pure strategy Nash equilibrium satisfying $\lambda f^{-1} = \tau$. The issue is that the pure strategy profile $f$ derived from $g$ may not actually be a Nash equilibrium, as it does not guarantee that $f(i) \in \cA(i)$. In fact, the existing literature on the closed graph property of large games primarily focuses on a common action space and rarely considers the action correspondence. As a result, the purification results in those works cannot be applied directly to our context.}

%-------------------------------------------------------------------------

\section{An application: equilibrium selection in large games}\label{sec-application}
In this subsection, we introduce an application of our main theorem. A large game may have multiple Nash equilibria, some of which may be unreliable. We propose an equilibrium selection criterion based on finite approximation. Since a large game is an idealization of (large) finite-player games, this idealization process may yield some Nash equilibria that exist only in the limit large game but not in the finite-player games. Such equilibria in the large game are unreliable and should be excluded. To illustrate this idea, we consider the following example.

\begin{exam}\label{exam-2}
Let $G^n$ be a large finite-player game with $n$ players and a common action set $A = \{a,b\}$. Players have a common payoff function $u$ given by
$$u(a, \tau) = \tau(a); \,\,\,\,\,\,\,\, u(b, \tau) = \tau(a) + 2\tau(b).$$
As $n$ tends to infinity, we obtain an idealized large game $G$ with a continuum of players. Players in $G^n$ may assume they are playing $G$ to simplify the analysis.

Clearly, the large game $G$ has two Nash equilibria: everyone choosing $a$ (denoted by $f_1$) and everybody choosing $b$ (denoted by $f_2$). However, once players realize that $G$ is not the actual game but that $G^n$ is the real game they are playing, they will see that action $a$ is strictly dominated by $b$ in $G^n$. Hence, any Nash equilibrium of $G^n$ does not assign a positive probability to action $a$. Therefore, $f_2$ is the only reliable Nash equilibrium of $G$, as it is a limit point of the Nash equilibria of finite-player games $\{G^n\}_{n \in \bZ_{+}}$.
\end{exam}

Equilibrium selection via finite approximation works well in Example~\ref{exam-2}; however, when generalizing this method to general large games, a natural question arises: will this selection process lead to an empty set in some large game? The result below, which is a direct corollary of Theorem~\ref{thm-closed graph property}, shows that a reliable Nash equilibrium that can be approximated by equilibria of large finite-player games always exists.

\begin{cro}\label{application}
Given a sequence of large finite-player games $\{G^n\}_{n \in \bZ_{+}}$ and the limit large game $G$, there exists a Nash equilibrium action distribution $\tau$ of $G$ that is a limit point of Nash equilibria of $\{G^n\}_{n \in \bZ_{+}}$. In particular, if $G$ has a unique Nash equilibrium action distribution $\tau^*$, then any  sequence of Nash equilibria from  $\{G^n\}_{n \in \bZ_{+}}$converges to $\tau^*$.  
\end{cro}

Corollary~\ref{application} serves as the theoretical foundation for our equilibrium selection criterion. Moreover, this selection process can be used iteratively to obtain more accurate predictions about the outcomes of large games. 

\begin{exam}\label{exam-3}
Let $G^n$ be a large finite-player game with $n$ players and a common action set $A = \{a,b,c,d\}$. Players have a common payoff function $u$ that is given by
	\begin{align*}
                 &u(a, \tau) = \tau(a)\\
			 &u(b, \tau) = \tau(a) + 2\tau(b) \\
                 &u(c, \tau) = \tfrac{\tau(a)}{2} + 2\tau(b) + 3\tau(c)\\
			 &u(d, \tau) = \tfrac{\tau(a)}{2} + \tau(b) + 3\tau(c) + 4\tau(d).
	\end{align*}
As $n$ tends to infinity, we obtain an idealized large game $G$ with a continuum of players. Clearly, the large game $G$ has multiple Nash equilibria, such as all players choosing the same action. Among these Nash equilibria, only one is reliable as a limit point of convergent Nash equilibria of finite-player games. In fact, once players realize they are playing a finite game $G^n$, they will observe that action $a$ is strictly dominated by action $b$ in $G^n$, meaning $a$ will not be chosen in a Nash equilibrium of $G^n$ ( i.e., $\tau(a) = 0$). Based on this observation, players will also see that $b$ is strictly dominated by $c$, leading to $\tau(b) = 0$ in any Nash equilibrium of $G^n$. Proceeding inductively, players will ultimately conclude that the only reliable Nash equilibrium is $f^*(i) \equiv d$. Thus, among all the Nash equilibria of $G$, only $f^*$ is a limit point of Nash equilibria of finite-player games.
\end{exam}

%-------------------------------------------------------------------------

\section{Related literature}\label{sec-literature}
Significant progress has been made in the study of games and economies with a continuum of agents since the foundational papers \cite{Aumann1964}, \cite{S1973}, and \cite{M1984}. Recent theoretical developments include contributions from \cite{Kalai2004},  \cite{KRSY2013}, \cite{QY2014}, \cite{Yu2014}, \cite{DK2015}, \cite{HSS2017}, \cite{KRYZ2017}, \cite{KS2018}, \cite{KRQS2020}, \cite{Hellwig2022},  \cite{CQSS2022}, \cite{ADKU2022a, ADKU2022b}, and \cite{Yang2022, Yang2023}. In addition to these theoretical advances, many papers have explored the convergence of equilibrium in various applied settings, such as matching and contests. Here, we will focus only on theoretical works that are directly related to this paper.

This paper is directly related to the literature on closed graph property of Nash equilibrium. This topic was explored in earlier works such as \cite{Green1984} and \cite{Housman1988}. Recent papers include \cite{KRSY2013}, \cite{QY2014}, \cite{QYZ2016},  \cite{HSS2017} and \cite{Wu2022}. The key distinction is the focus on the convergent sequence of finite-player games: whereas these earlier papers concentrated on pure strategy Nash equilibria, our paper extends the setting by considering randomized strategy Nash equilibria.

Technically, this paper also relates to existing research on randomized strategy profile/Nash equilibrium in large finite-player games. Relevant papers include \cite{Kalai2004}, \cite{DK2015}, and \cite{KS2018}. However, the research questions in these papers differ from ours: the first two papers studied the relationship between a randomized strategy Nash equilibrium and its possible realizations. Based on some assumptions of equicontinuity or Lipschitz continuity of payoff functions, they proved that a randomized strategy Bayes-Nash equilibrium is approximately ex post Nash in large finite-player games. \cite{KS2018} examined the equilibrium property of an imagined continuum symmetric equilibrium in a repeated game with incomplete information, where payoff functions also satisfy some Lipschitz continuity condition.
%\footnote{See \cite{PY2014} for the existence of randomized strategy Nash equilibrium in discontinuous games; also see \cite{LY2022} and \cite{LSY2020} for randomized strategy implementation under ambiguity.} %-------------------------------------------------------------------------

\section{Appendix}\label{sec-appendix}
%-------------------------------------------------------------------------
\subsection{Proofs of Claims in Subsection~\ref{subsec-motivation}}\label{subsec-proof of claim}
\begin{proof} [Proof of Claim~\ref{claim-no PSNE}]
	Suppose that the game $G^n$ possesses a pure strategy Nash equilibrium  $f^n \colon I^n \to A$. Let $\tau^n(a)$ denote the proportion of players choosing the action $a \in \{R, S, P\}$ under the Nash equilibrium $f^n$, and assume that $\tau^n(R)= x,\tau^n(S)= y, \tau^n(P)= z,$  obviously $x, y, z \in [0,1]$ and $x + y + z = 1$. Without loss of generality, let $ x = \text{min}\{x, y, z\}$  and hence we only need to discuss the cases $x \le y \le z$ and $x \le z \le y$.

	For the case $x \le y \le z$, we divide the discussion into the following three parts.
	\begin{itemize}
		\item[(i)] $x = y = 0, z =1$.\\
		In this case the payoff of the player who chooses $P$ is $0$. However, if the player choosing $P$ deviates to $S$, then her payoff will be $1 - \frac{1}{n} > 0$. Thus, players choosing action $P$ have an incentive to deviate and hence $f^n$ cannot be a pure strategy Nash equilibrium.
		\item[(ii)] $x = 0, 0 < y \le z$.\\
		In this case the payoff of the player who chooses $P$ is $-y < 0$. However, if the player choosing $P$ deviates to $S$, then her payoff will be $z - \frac{1}{n} \geq 0$. Thus, players choosing action $P$ have an incentive to deviate and hence $f^n$ cannot be a pure strategy Nash equilibrium.
		\item[(iii)] $0 < x \le y \le z$.\\
		In this case the payoff of the player who chooses $R$ is $y - z \le 0$. However, if the player choosing $R$ deviates to $S$, then her payoff will be $z - x + \frac{1}{n} > 0$.  Thus, players choosing action $R$ have an incentive to deviate and hence $f^n$ cannot be a pure strategy Nash equilibrium.
	\end{itemize}
	
For the case $ x \le z \le y$, we can  similarly divide the discussion into three parts, and we can see that there always exist some players who have an incentive to unilaterally deviate. In conclusion, there does not exist any pure strategy Nash equilibrium in this game when $n \ge 2$. 

\end{proof}

\begin{proof} [Proof of Claim~\ref{claim-large symmetric}]
Clearly, we have $u(R, \tau^*) = u(S, \tau^*) = u(P, \tau^*) = \tfrac{1}{3} - \tfrac{1}{3} = 0$. Since every single player is negligible in the large game $G$, we can see that all the players will have the same payoff by choosing any one of the actions in $\{R,S,P\}$. Hence given the societal summary $\tau^*$, no player has an incentive to deviate and $g$ is a  Nash equilibrium of $G$.

Then we prove the uniqueness of the equilibrium action distribution. Assume $\tau$ is an equilibrium action distribution of $G$, and $\tau(R)=x, \tau(S)=y, \tau(P)=z$.
	 It is clear that $x, y, z \in [0,1]$ and $x + y + z = 1$. Without loss of generality, let $ x = \text{min}\{x, y, z\}$  and hence we only need to discuss two cases, i.e. $x \le y \le z$ and $x \le z \le y$. \\	
	For the case $x \le y \le z$,  we divide the discussion into the following two parts.
	\begin{itemize}
		\item[(i)] $x < y $ or $y < z$.\\
		In this case the portion of players choosing action $P$ with positive probability is non-negligible. Since action $P$ is strictly dominated by action $S$ for each player ($x - y  \le 0 < z - x $), those players who choose action $P$ with positive probability do not play the optimal strategy. Thus, $\tau$ cannot be an equilibrium  action distribution.
		\item[(ii)] $x = y  = z = \frac{1}{3}$.\\
		It is easy to see that $\tau = (\frac{1}{3}, \frac{1}{3}, \frac{1}{3})$ is an equilibrium action distribution.
	\end{itemize}
    For the case $ x \le z \le y$, we can  similarly divide the discussion into two parts, and conclude that $\tau = (\frac{1}{3}, \frac{1}{3}, \frac{1}{3})$ is the unique equilibrium action distribution. 
\end{proof}

\begin{proof} [Proof of Claim~\ref{claim-finite symmetric}]
To show that $g^n (i) \equiv \frac{1}{3}\delta_{R} + \frac{1}{3}\delta_{S} + \frac{1}{3}\delta_{P}$ is a  Nash equilibrium of $G^n$, we need to prove that every player $i$ is indifferent between the actions in $\{R,S,P\}$, given other players follow the strategy profile $g^n_{-i}$. We first calculate player $i$'s (expected) payoff if she deviates to action $R$ as follows:
	\begin{align*}
                  &u(R,g_{-i}^n) = \int_{\prod\limits_{j =1}^{n-1} A} u \Bigl( R, \frac{1}{n} \bigl(\sum\limits_{j \in I^n \backslash \{i\}} \delta_{a_j} +\delta_{R} \bigr)\Bigr) \tens{j \in I^n \backslash \{i\}} g^n(j, \rmd a_j) \\
			 &= \int_{\prod\limits_{j =1}^{n-1} A}      \Bigl(\frac{1}{n} \sum\limits_{j \in I^n \backslash \{i\}} \delta_{a_j}  \Bigr)(S)   -  \Bigl( \frac{1}{n}\sum\limits_{j \in I^n \backslash \{i\}} \delta_{a_j} \Bigr)(P)  \tens{j \in I^n \backslash \{i\}} g^n(j, \rmd a_j ) \\
                &=  \int_{\prod\limits_{j =1}^{n-1} A}      \Bigl( \frac{1}{n}\sum\limits_{j \in I^n \backslash \{i\}} \delta_{a_j}  \Bigr)(S) \tens{j \in I^n \backslash \{i\}} g^n(j, \rmd a_j ) - \int_{\prod\limits_{j =1}^{n-1} A} \Bigl( \frac{1}{n}\sum\limits_{j \in I^n \backslash \{i\}} \delta_{a_j} \Bigr)(P)  \tens{j \in I^n \backslash \{i\}} g^n(j, \rmd a_j ) \\
			&= 0,
	\end{align*}
where the second equality follows from the fact that $u(R, \tau) = \tau(S) - \tau(P)$ and $\delta_R(S) = \delta_R(P) = 0$, and the last equality holds as players choose $S$ and $P$ with equal probability. By the same argument, we have that $u(S,g_{-i}^n) = u(P,g_{-i}^n) =0$. Therefore, player $i$ has no incentive to deviate and $g$ is a Nash equilibrium of $G$.
    \end{proof}

\begin{proof} [Proof of Claim~\ref{claim-invest}]
We prove by contradiction. Assume that $\tau^* = (\tfrac{1}{2}, \tfrac{1}{2})$ is a limit point of $\{g^n\}_{n \in \bZ_{+}}$, where $g^n$ is a pure strategy Nash equilibrium of $G^n$. For each $n \ge 2$, let $\tau^n$ denote the action distribution induced from $g^n$. Without loss of generality, we assume that $\lim_{n \to \infty} \tau^n = \tau^*$. Thus, for sufficiently large $n$, we have $\min\{\tau^n(a_1), \tau^n(a_2)\} > \tfrac{1}{n}$. If $\tau^n(a_1) \le \tau^n(a_2)$, then a player choosing $a_1$ would have an incentive to deviate to $a_2$ due to the strict monotonicity of the function $h$:
$$u(a_1, \tau^n) = h\bigl(\tau^n(a_1)\bigr) \le h\bigl(\tau^n(a_2)\bigr) < h\bigl(\tau^n(a_2) + \tfrac{1}{n}\bigr)  = u(a_2, \widetilde\tau^n),$$
where $\widetilde\tau^n = \tau^n - \tfrac{1}{n}\delta_{a_1} + \tfrac{1}{n}\delta_{a_2}$ represents the action distribution after the player's deviation. Hence, $g^n$ is not a Nash equilibrium, leading to a contradiction. On the other hand, if $\tau^n(a_1) > \tau^n(a_2)$, we can similarly conclude that players choosing $a_2$ would deviate, which also leads to a contradiction. 
\end{proof}
%%%%%%%%%%%%%%%%%%%%%%%%%%%%%%%%%%%%%%%%%%%%%%%%%%%%%%%%%%%%%%%%%%%%%%%%%%%%%%%%%%%%%%%%%%%%%%%%%%%%%%%%%%%%%%%%%%%%%%%%%%

\subsection{Proof of the main theorem}\label{subsec-proof main results}
Our proof of Theorem~\ref{thm-closed graph property} is divided into two parts. In Subsection~\ref{subsubsec:lem}, we establish two lemmas as the technical preparation. Using these lemmas, we complete the proof  of Theorem~\ref{thm-closed graph property} in Subsection~\ref{subsubsec:main}.

\subsubsection{Lemmas}\label{subsubsec:lem}
Let $d_A$ be the metric on the action space $A$ such that the Borel $\sigma$-algebra induced from $d_A$ is $\cB(A)$. Given a sequence of randomized strategy profiles $\{g^n\}_{n \in \bZ_+} \colon I \to \cM(A)$, let $\{x^n_i\}_{ i \in I^n, n \in \bZ_{+}}$  be a sequence of random variables mapping from a probability space $(\Omega, \Sigma, \mathbb{P})$ to $A$ such that (i) for each player $i \in I^n$ and $n \in \mathbb{Z}_{+}$, the distribution induced from $x^n_i$ is $g^n(i)$;\footnote{That is, $\bP(x^n_i)^{-1} = g^n(i)$ for each $i \in I^n$ and $n \in \bZ_+$.} and (ii) for each $n \in \mathbb{Z}_{+}$, the random variables $\{x^n_i\}_{i \in I^n}$  are pairwise independent.

Lemma~\ref{lem-norm}, taken from \cite{CWX2024}, shows that given a sequence of strategy profiles of finite-player games, the distance between the realized societal summary and the ex ante society summary converges to 0 in probability.

\begin{lem}\label{lem-norm}
	Let $\{G^{n} \}_{ n \in \mathbb{Z}_{+}}$ be a sequence of  finite-player games, and $g^n$ be a randomized strategy profile of $G^n$ for each $n \in \mathbb{Z}_{+}$. For each $\omega \in \Omega$, let $s(x^n)(\omega) = {\sum_{i \in I^n}{{\delta _{ x^n_i(\omega) }\lambda^n(j)}} }$ be a realized societal summary of the randomized strategy profile $g^n$, hence $s(x^n)$ can be viewed as a random variable mapping from  $(\Omega, \Sigma, \bP) $ to $\cM(A)$. Then we have
	$$\rho\bigl(s(x^n), s(g^n)\bigr) \to 0 \text{   and   } \beta \bigl( s(x^n), s(g^n)\bigr) \to 0\text{  in probability,}$$
where $\rho$ denotes the Prohorov metric on $\cM(A)$,\footnote{Definition of the Prohorov metric $\rho$: for all $\tau,\widetilde \tau \in \mathcal{M}(A)$, we have
	\[\rho( \tau, \widetilde \tau )=\inf \bigl\{\varepsilon>0 \colon \tau(B) \leqslant \varepsilon + \widetilde \tau (B^{\varepsilon}) , \widetilde \tau(B) \leqslant \varepsilon + \tau (B^{\varepsilon}) \text { for all  } B \in \cB(A) \bigr\},\] where
	${B^\varepsilon } = \bigl\{ a \in A \colon d_A(a,b) < \varepsilon {\text{ for some }} b \in B \bigr\} .$} and $\beta$ represents the dual-bounded-Lipschitz metric on $\mathcal{M}(A)$.\footnote{Definition of the dual-bounded-Lipschitz metric $\beta$: for all   $\tau,\widetilde \tau \in \mathcal{M}(A)$, we have
	\begin{equation*}
		\beta(\tau, \widetilde \tau) = \| \tau-\widetilde \tau\|_{{B L}}^{*}=\sup \Bigl\{\Bigl|\int_{A} h {\mathrm{d}}(\tau- \widetilde \tau)\Bigr| \colon\|h\|_{{B} {L}} \leqslant 1 \Bigr\},
	\end{equation*}
	where $h$ is bounded continuous on $A$,
	$\|h\|_{\infty} = \sup\limits_{a \in A}|h(a)|$, $\|h\|_{\mathrm{L}}=\sup\limits_{a \neq b, a,b \in A } \dfrac{|h(a)-h(b)|}{d_A(a, b)}$,  and $\|h\|_{{BL}}  = \|h\|_{\infty}+\|h\|_{{L}}.$}
\end{lem}

Notice that the Prohorov metric $\rho$ and the dual-bounded-Lipschitz metric $\beta$ are equivalent metrics on $\cM(A)$(\citet[Theorem 8.3.2]{Bo2007}). The detailed proof of Lemma~\ref{lem-norm} is in \citet[Subsection 6.1]{CWX2024}. Based on Lemma~\ref{lem-norm}, we prove an inequality in Lemma~\ref{lem-payoff gap} below showing that a Nash equilibrium $g^n$ of $G^n$ is approximately Nash when the aggregated action distribution is exactly equal to the society summary $s(g^n) = \sum_{j \in I^n} \lambda^n(j)g^n(j)$. Note that for any randomized strategy profile $g^n$,  $u^n_i \Bigl( {{g^n}(i),\sum\limits_{j \in I^n} \lambda^n(j){{g^n}} (j)} \Bigr)  =  \int_{A_i^n}   u^n_i \Bigl( {a,\sum\limits_{j \in I^n} \lambda^n(j){{g^n}} (j)} \Bigr) {g^n}(i, \rmd a)$.

\begin{lem}\label{lem-payoff gap}
	Let $\{G^{n} \}_{ n \in \mathbb{Z}_{+}}$ be a sequence of  finite-player games converges weakly to a large game $G$. Given any $\gamma > 0$, there exists a  sequence of sets  $ S^n \subseteq I^n$  such that
		\begin{itemize}
			\item[(i)]   $ \lambda^{n}(S^n) > 1- \gamma$ for all $n \in \mathbb{Z}_+$, and $\{ u_i^n  \}_{ i \in  S^n, n \in \mathbb{Z}_+}$ are equicontinuous and uniformly bounded by a constant $M_\gamma$.
			\item[(ii)]  For any $\varepsilon > 0$ and any sequence of Nash equilibira  $\{g^n\}_{n \in \mathbb{Z}_+}$ (each $g^n$ is a Nash equilibrium of $G^n$ ), there exists $N_{\varepsilon} $ such that
			\[
			u^n_i \Bigl( {{g^n}(i),\sum\limits_{j \in I^n} \lambda^n(j){{g^n}} (j)} \Bigr)
			\ge
			u^n_i \Bigl( a, \lambda^n(i) \delta_{a} + \sum\limits_{j \in I^n \backslash \{i\}} \lambda^n(j)g^n(j) \Bigr) - \frac{\varepsilon }{2}
			\]
			for all $n \geq N_{\varepsilon}$,  $ i \in S^n ,a \in A_i^{n}$.
	\end{itemize}
\end{lem}

\begin{proof}[Proof of Lemma~\ref{lem-payoff gap}]
Our proof consists of three major steps. In step 1, we construct the sequence of sets $\{S^n\}_{n \in \mathbb{Z}_+ }$ that satisfies the first requirement. Then we establish an estimation of the difference between  $u_i^n( \mu, g^n_{-i})$ and  $u^n_i \bigl( \mu , \lambda^n(i) \mu + \sum_{j \in I^n \backslash \{i\} } \lambda^n(j) g^n (j) \bigr)$ for all $i \in  S^n$, $\mu \in \cM(A)$.  Finally in step 3, we prove the inequality in part (ii) of the Lemma~\ref{lem-payoff gap} by using part (i) and the estimation result in proved step 2.

\bigskip
	
\noindent\textbf{Step 1.} 	Notice that the game $G^n$ (resp. $G$) is a mapping from the player set $I^n$ (resp. $I$) to $\cC_A \times \cU_A$. Thus $G^n(i) = \bigl(G^n_1(i), G^n_2(i)\bigr)$ and $G(i) = \bigl(G_1(i), G_2(i)\bigr)$. Let $\mathcal{W}^{n} = \lambda^n  (G_{2}^n)^{-1} $ and $\mathcal{W} = \lambda  (G_{2})^{-1}$. Clearly, both $\mathcal{W}^{n}$ and $\mathcal{W}$ are probability measures on $\cU_A$. Since $A \times \cM(A)$ is a compact metric space, the space of bounded and  continuous functions $\cU_A$ on $A \times \cM(A)$ is a Polish space. By using the Prohorov theorem (\citet[Theorem 5.2]{B1999}), we know that $\{ \mathcal{W}^n \}_{n \in \bZ_{+}}$ is tight, which means that for any $ \gamma > 0$, there exists a compact set $K_{\gamma} \subset \cU_A$ such that $\mathcal{W}^{n}(K_{\gamma}) > 1- \gamma   $ for all $n \in \bZ_{+} $.

Since $K_{\gamma}$ is a compact set that consists of bounded and continuous functions, the Arzelà-Ascoli theorem (\citet[Theorem 45.4]{Mu2000}) implies that all the functions in $K_{\gamma}$ are equicontinuous and uniformly bounded. Let $M_\gamma$ denote a bound of all the functions in $K_{\gamma}$, and $ S^n = \{ i \in I^n |  u_i^n \in K_{\gamma}  \}$ for all $n \in \mathbb{Z}_+$. Obviously, we have $\lambda^n( S^n ) = \mathcal{W}^{n}(K_{\gamma})> 1  -\gamma$.

\bigskip

\noindent\textbf{Step 2.} Let $\{x^n_i\}_{ i \in I^n, n \in \bZ_{+}}$  be a sequence of random variables mapping from a probability space $(\Omega, \Sigma, \mathbb{P})$ to the action space $A$ such that each induced action distribution $\bP (x^n_i)^{-1} = g^n(i)$. Fix a player $i$, pick arbitrary a randomized strategy $\mu \in \cM(A)$ and let $ x_{\mu}$ be a random variable maps from $\Omega$ to $A$ that induces the distribution $\mu$. Then we can rewrite player i's expected payoff as follows:
$$u_i^n( \mu, g^n_{-i}) = \bE \Bigl[   u_i^n \Bigl(   x_{\mu} , \lambda^n(i)\delta_{x_{\mu}} +  \sum_{j \in I^n \backslash \{i\} } \lambda^n(j)\delta_{x^n_j}       \Bigr) \Bigr],$$
and similarly:
$$u^n_i \Bigl( \mu , \lambda^n(i) \mu +  \sum_{j \in I^n \backslash \{i\}} \lambda^n(j)g^n (j) \Bigr)
=
\bE \Bigl[   u_i^n\Bigl(   x_{\mu},  \lambda^n(i) \mu +  \sum_{j \in I^n \backslash \{i\}} \lambda^n(j)g^n (j)    \Bigr)\Bigr].$$

Hereafter, we focus on the payoff functions in $K_{\gamma}$. By equicontinuity, we know that for any $\varepsilon > 0$, there exists $\eta>0$ such that for any $ \tau, \widetilde{\tau} \in \cM(A)$  with $\rho( \tau, \widetilde \tau) \le \eta$, any $u \in K_{\gamma}$, and any $a \in A$, we must have
\begin{equation}\label{equa-payoff equi}
	\begin{split}
	|  u ( a , \tau )   - u ( a , \widetilde \tau)              |
		\le
		\frac{\varepsilon}{4(2M_\gamma+1)}.
	\end{split}
\end{equation}
Let $\cS (x_{\mu}, x^n_{-i})   (\omega) =  \lambda^n(i) \delta_{x_{\mu}(\omega)} +  \sum_{j \in I^n \backslash \{i\}   } \lambda^{n}(j)  \delta _{
	x^n_j(\omega)}  $, and $\cS (\mu, g^n_{-i})    =  \lambda^n(i) \mu +  \sum_{j \in I^n \backslash \{i\}  } \\ \lambda^{n}(j) g^n(j)  $, for all $\omega \in \Omega$, $\mu \in \cM(A)$.
The triangle inequality implies that
\begin{align}\label{align-metric triangle}
 \rho\bigl( \cS  (\mu, g^n_{-i})   ,  \cS  (x_{\mu}, x^n_{-i})   (\omega)  \bigr)  &  \le  \rho \bigl(  s(g^n) ,  \cS  (\mu, g^n_{-i})\bigr) 
 + \rho\bigl(   s(x^n)(\omega),  s(g^n)\bigr)  \notag \\ 
  & 	+
\rho\bigl(  s(x^n)(\omega)  ,  \cS  (x_{\mu}, x^n_{-i})   (\omega)\bigr).
\end{align}
By the definition of $\rho$, we know that $\rho\bigl( s(x^n)(\omega)   ,  \cS (x_{\mu}, x^n_{-i})   (\omega)\bigr) \le  \sup_{j \in I^{n}} \lambda^{n}(j)$ for any $\omega \in \Omega$, $\mu \in \cM(A)$, $i \in I^n$. Since $\sup_{j \in I^{n}} \lambda^{n}(j) \rightarrow 0 $, there exists $N_1 \in \mathbb{Z}_+$ such that for any $ n \ge N_1$, we have $\sup_{ j\in I^{n}} \lambda^n(j) < \frac{\eta}{4} $. Hence for any $n \ge N_1$, $i \in I^n$, $\mu \in \cM(A)$, $\omega \in \Omega$,
\begin{align}\label{equa-metric gap 1}
	\rho\bigl( s(x^n)(\omega)   ,   \cS(x_{\mu}, x^n_{-i})   (\omega)\bigr)  < \frac{\eta}{4}.
\end{align}
By the same argument as above, we can see that for any $n \ge N_1$, $i \in I^n$, $\mu \in \cM(A)$,
\begin{align}\label{equa-metric gap 2}
	\rho\bigl(s(g^n) ,    \cS (\mu, g^n_{-i})\bigr)   < \frac{\eta}{4}.
\end{align}

Let  $\Omega_1^{(\frac{\eta}{2},n)} = \bigl\{   \omega \in \Omega    \,|\,    \rho \bigl(  s(x^n)(\omega),  s(g^n)  \bigr) 	 < \frac{\eta}{2}    \} $ and $\Omega_2^{(\frac{\eta}{2},n)} = \Omega \backslash \Omega_1^{(\frac{\eta}{2},n)}$.
By Lemma~\ref{lem-norm}, for any $\varepsilon > 0 $, there exists $N_{\varepsilon} \ge N_1$ such that for any $n \ge N_{\varepsilon}$, we have $\mathbb P \Bigl( \Omega_2^{(\frac{\eta}{2},n)} \Bigr)   \le   \frac{\varepsilon }{4(2M_\gamma + 1)}$. Let
$$ H_1^{(\frac{\eta}{2}, n)} =
\left| \bE \Bigl[ \Bigl( u_i^n \bigl(   x_{\mu} ,    \cS (\mu, g^n_{-i}) \bigr)
-
u^n_i \bigl( x_{\mu} ,  \cS (x_{\mu},x^n_{-i}) \bigr) \Bigr) \delta_{ \Omega_1^{(\frac{\eta}{2},n)}  }  \Bigr]
\right|,$$
and
$$ H_2^{(\frac{\eta}{2}, n)}
=
\left| \bE \Bigl[ \Bigl( u_i^n \bigl(   x_{\mu} ,    \cS (\mu, g^n_{-i})  \bigr)
-
u^n_i \bigl( x_{\mu} ,  \cS (x_{\mu},x^n_{-i}) \bigr) \Bigr) \delta_{ \Omega_2^{(\frac{\eta}{2},n)}  }  \Bigr]
\right|.$$
By using the triangle inequality, we have
\begin{equation*}
	\begin{split}
		\left| \bE \Bigl[  u_i^n \bigl(   x_{\mu} ,    \cS (\mu, g^n_{-i})    \bigr)
		-
		u^n_i \bigl( x_{\mu} ,  \cS( x_{\mu},x^n_{-i}) \bigr)   \Bigr]
		\right|
		\leq H_1^{(\frac{\eta}{2}, n)} + H_2^{(\frac{\eta}{2}, n)}.
	\end{split}
\end{equation*}

Then we estimate $H_1^{(\frac{\eta}{2}, n)}$ and $H_2^{(\frac{\eta}{2}, n)}$ separately.  For each $n \ge N_{\varepsilon}$ and player $i \in S^n$, we have $u_i^n \in K_{\gamma}$. By the definition of event $\Omega_1^{(\frac{\eta}{2}, n)}$ and the Inequalities~(\ref{equa-payoff equi}), (\ref{align-metric triangle}), (\ref{equa-metric gap 1}), and (\ref{equa-metric gap 2}), we can see that
	\begin{align}\label{equa-payoff-H1}
		H_1^{(\frac{\eta}{2}, n)} \le \frac{\varepsilon }{4(2M_\gamma + 1)}.
	\end{align}
Since $u_i^n$ is bounded by $M_\gamma$, we have
	\begin{align}\label{equa-payoff-H2}
		H_2^{(\frac{\eta}{2}, n)} \le  2M_\gamma\frac{\varepsilon }{4(2M_\gamma + 1)}.
	\end{align}
Combining Inequalities (\ref{equa-payoff-H1}) and (\ref{equa-payoff-H2}) we can see that
$$\left|	u_i^n( \mu, g^n_{-i}) - u^n_i \Bigl( \mu , \lambda^n(i) \mu + \sum\limits_{j \in I^n \backslash \{i\} } \lambda^n(j) g^n (j) \Bigr) \right| \le \frac{\varepsilon}{4}.$$

\bigskip

\noindent\textbf{Step 3.} Now we know that for any $\varepsilon>0,$ there exists $N_{\varepsilon} $ such that  for all $n \ge N_{\varepsilon}, i \in S^n$, $a \in A$, we have
\begin{align}\label{equa-payoff dif 1}
	\Bigl|	 u^n_i \Bigl( g^n(i),\sum\limits_{j \in I^n} \lambda^n(j) g^n (j) \Bigr) - u_i^n(g^n)  \Bigr| \le \frac{\varepsilon}{4},
\end{align}
and
\begin{align}\label{equa-payoff dif 2}
	\Bigl|	 u^n_i \Bigl( a, \lambda^n(i) \delta_{a} + \sum\limits_{j \in I^n \backslash \{i\}} \lambda^n(j) g^n (j) \Bigr) - u_i^n(a,g^n_{-i})  \Bigr| \le \frac{\varepsilon}{4}.
\end{align}

Thus for any $n \geq N_{\varepsilon}$, $i \in S^n$, $a \in A_i^n$, we have
\begin{align*}%\label{equa-lem3}
	u_i^n \Bigl(g^n(i),  \sum\limits_{j \in I^n} \lambda^n(j) g^n(j)\Bigr) 
	&\ge
	u_i^n (g^n) - \frac{\varepsilon }{4} \\
	&\ge u_i^n(a,g^n_{-i})  - \frac{\varepsilon }{4} \\
	&\ge   u^n_i\Bigl( a, \lambda^n(i) \delta_{a} + \sum\limits_{j \in I^n \backslash \{i\}} \lambda^n(j) g^n (j) \Bigr) - \frac{\varepsilon }{2},
\end{align*}
where the first inequality follows from (\ref{equa-payoff dif 1}), the second inequality follows  from the definition of Nash equilibrium and the third inequality follows from (\ref{equa-payoff dif 2}). 
\end{proof}

\subsubsection{Proof of Theorem~\ref{thm-closed graph property}}\label{subsubsec:main}
We first construct a strategy profile $g \colon I \to \cM(A)$ such that $s(g) = \tau^*$. Recall that $G^n(i) = \bigl( \cA^n(i), u^n_i\bigr) = (A^n_i, u^n_i)$ for each $n \in \bZ_{+}$ and each player $i \in I^n$. $\{g^n\}_{n \in \bZ_{+}}$ is a sequence of Nash equilibria of $\{ G^n \}_{n \in \bZ_{+}}$
such that the societal summaries $ \{ \int_{ I^n}  g^n \mathrm{d} \lambda^n \}_{ n \in \bZ_{+}}$ converges weakly to some $\tau^* \in \mathcal{M}(A)$. For each $n \in \mathbb{Z}_{+}$, let $\nu^n=\int_{{I^n}} \delta_{G^n(i)}  \otimes {g^n}(i){\rm{d}}{\lambda^n}(i)$ be a joint measure on $\cC_A \times \cU_A \times A$. By compactness, there exists a subsequence of $\{\nu^n\}_{n \in \bZ_+}$ that weakly converges to $\nu$, and without loss of generality, we simply assume that the whole sequence $\{\nu^n\}_{n \in \bZ_+}$ converges to $\nu$.

Since the sequence of games $\{G^n\}_{n \in \bZ_+}$ weakly converges to $G$, we know that the sequence $\nu^n|_{\cC_A \times \cU_A} = \lambda (G^n)^{-1}$ weakly converges to $\nu|_{\cC_A \times \cU_A} = \lambda G^{-1}$. As $\cC_A \times \cU_A$ and $A$ are Polish spaces, there exits a family of Borel probability measures $\{   \widetilde{\nu}( B, u, \cdot    ) \}_{ (B, u) \in \cC_A \times \cU_A}$ in $\cM(A)$, which is the disintegration of $\nu $ with respect to $\lambda G^{-1}$. Let $g$ be a measurable function from $I$ to $\cM(A)$ such that
$$g(i,Q) = \widetilde{\nu} ( G(i), Q  )$$
for any $i \in I$, $Q \in \cB(A)$. Hence we have $\nu = \int_{ I  } \delta_{G(i)} \otimes g(i) \rmd \lambda (i)$.
Since the sequence $\{s(g^n)= \int_{ I^n} g^{n} \mathrm{d} \lambda^n \}_{n \in \bZ_{+}}$ weakly converges to $\tau^*$, and $\{\nu^n\}_{n \in \bZ_+}$ weakly converges to $\nu$, we conclude that $s(g) = \int_I g \rmd \lambda = \tau^*$. Thus it suffices to show that $g$ is a  Nash equilibrium of $G$. We divide the remaining proof into three steps. In step 1, we show that $\mathrm{supp\,\,} g(i) \subseteq A_i$ for all $i \in I$,  where $\mathrm{supp\,\,} g(i)$ is the smallest closed set $B \subseteq A$ such that $g(i, B) = 1$. In step 2, we construct an auxiliary function $\Psi$ that plays an essential role in our proof. Finally in step 3, we show that $g$ is a Nash equilibrium of the large game $G$.

\bigskip

\noindent\textbf{Step 1.} Let $\cS( \cA^n , g^n  ) = \int_{I^n} \delta_{ \cA^n(i)} \otimes g^n(i) \rmd \lambda^n(i)$ for each $n \in \bZ_{+}$, and $ \cS(\cA,g) = \int_{I} \delta_{ \cA(i)} \otimes g(i) \rmd \lambda(j)$. Since $\nu^n$ converges weakly to $\nu$, $\nu^n|_{\cC_A \times A}$ also converges weakly to $\nu|_{\cC_A \times A}$, which implies that  $\cS( \cA^n , g^n  )$ converges weakly to $\cS( \cA , g)$.
Let $Z = \{ (B,b) | B \in \cC_A, b \in B \}$. Clearly, $Z$ is a closed subset and
$\cS(\cA^n, g^n)(Z) = 1 $. By the weak convergence, we have $\cS(\cA, g)(Z) = 1.$ By the definition of $Z$, we can see that $\mathrm{supp\,\,} g(i)  \subseteq A_i$ for $\lambda$-almost all $i \in I$.

\bigskip

\noindent\textbf{Step 2.} Notice that $\nu^n|_{ A}=s(g^n)$ for each $n \in \mathbb{Z}_{+}$.
Since $g^n $ is a Nash equilibrium of $G^n$, by Lemma~\ref{lem-payoff gap}, fix $\gamma > 0$ , for any $\varepsilon > 0$, there exists $N_{\varepsilon} \in \mathbb{Z}_{+}$ such that for any $n \geq N_{\varepsilon}$,
\[
u^n_i \bigl( g^n(i), s(g^n) \bigr)
\ge
u^n_i \bigl(a, s(a,g^n_{-i})\bigr) - \frac{\varepsilon }{2}
\]
for all $ i \in S^n$  and $a \in A_i^{n}$. Note that $S^n$ is a subset of $I^n$ with $ \lambda^{n}(S^n) > 1- \gamma,$ and
$s (a, g^n_{-i})   =  \lambda^n(i) \delta_{ a} +  \sum_{j \in I^n \backslash \{i\}  } \lambda^{n}(j)  g^n(j)  $.
Since $\sup_{j \in I^n} \lambda^n(j) \to 0$ and $\{s(  g^n   ) \}_{ n \in \mathbb{Z}_{+}  }$ converges weakly to $\tau$, $\{s (a, g^n_{-i})   \}_{ n \in \mathbb{Z}_{+}  }$ also converges weakly to $\tau$.

Let 
$$\cS \bigl(G^n,g^n,s(g^n)\bigr)  = \int_{I^n} \delta _{G^n(i)} \otimes g^n(i) \otimes \delta_{s(g^n)}  {\rm{d}}{{\lambda ^n}(i)}$$
for each $n \in \bZ_{+}$ and 
$$\cS \bigl(G,g,s(g)\bigr)  = \int_{I}  \delta _{G(i)} \otimes g(i) \otimes \delta_{s(g)}   \rmd \lambda (i).$$
Then \citet[Theorem 3.9]{B1999} implies that  $\{\cS (G^n,g^n,s(g^n)) \}_{ n \in \mathbb{Z}_{+}  }$ converges weakly to $ \cS(G,g,s(g))$.

Consider a mapping $\Psi \colon \cC_A \times \cU_A \times \cM(A) \times \cM(A) \to \bR$ defined as follows:
$$\Psi   (\widetilde{B}, \widetilde u, \widetilde \mu, \widetilde \tau) = \min_{a' \in \widetilde{B}} \left\{\int_A {\widetilde u(a, \widetilde\tau)} \widetilde \mu(\rmd a) - \widetilde u(a', \widetilde \tau)\right\}$$ 
for any $(\widetilde{B}, \widetilde u, \widetilde \mu, \widetilde \tau )  \in  \cC_A \times \cU_A \times \cM(A) \times \cM(A) $. Below we show that $\Psi$ is a continuous function on $\cC_A \times \cU_A \times \cM(A) \times \cM(A)$.

In fact, $\Psi$ can be rewritten as 
$$\Psi   ( \widetilde B, \widetilde u, \widetilde \mu, \widetilde \tau ) = \int_A {\widetilde u(a,\widetilde \tau )} \widetilde \mu ( \rmd a) - \max_{  \widetilde a \in \widetilde B } \widetilde u (\widetilde a, \widetilde \tau ) $$ 
for any $(\widetilde B, \widetilde u, \widetilde \mu, \widetilde \tau )  \in \cC_A \times  \cU_A \times \cM(A) \times \cM(A) $. Let $\psi \colon  \cU_A \times \cM(A) \times \cM(A) \to \bR$ be
$$\psi   ( \widetilde u, \widetilde \mu, \widetilde \tau ) = \int_A {\widetilde u(a,\widetilde \tau )} \widetilde \mu ( \rmd a)$$ 
for any $( \widetilde u, \widetilde \mu, \widetilde \tau )  \in   \cU_A \times \cM(A) \times \cM(A) $. We first verify that $\psi$ is continuous.
For any $(u', \mu', \tau'),(\widetilde u, \widetilde \mu, \widetilde \tau) \in \cU_A \times \cM(A) \times \cM(A) $,
\begin{align*}
	\psi(u',\mu',\tau') - \psi(\widetilde u,\widetilde \mu,\widetilde \tau) & = \int_A { u'(a, \tau' )}  \mu'( \rmd a) - \int_A {\widetilde u(a,\widetilde \tau )} \widetilde \mu( \rmd a) \\
	& = \int_A \bigl(u'(a, \tau' ) - {\widetilde u(a,\tau' )}\bigr)  \mu'( \rmd a) \tag{i} \\
	&  + \int_A \bigl(\widetilde u(a, \tau' ) - {\widetilde u(a, \widetilde \tau )}\bigr)  \mu'( \rmd a)  \tag{ii} \\
	& +  \int_A \widetilde u(a,\widetilde \tau )  (\mu'  - \widetilde \mu)( \rmd a).   \tag{iii}
\end{align*}
Part (i) tends to $0$ as $ ||u' - \widetilde u||_{\infty} \to 0  $. Since $A \times \cM(A)$ is compact, $\widetilde u$ is uniformly continuous, and part (ii) tends to $0$ when $\rho(\tau' , \widetilde \tau ) \to 0$. Finally, as $ \widetilde u(\cdot, \widetilde \tau) $ is a bounded continuous function on $A$, we know that part (iii) tends to $0$ as $\rho( \mu' , \widetilde \mu) \to 0$. Thus, $\psi$ is continuous. Let $\phi \colon \cU_A \times \cC_A \times \cM(A) \to \bR$ be 
$$\phi( \widetilde u, \widetilde B, \widetilde \tau   ) =   \max_{  \widetilde a \in \widetilde B } \widetilde u (\widetilde a, \widetilde \tau )$$ 
for any  $(\widetilde u, \widetilde B, \widetilde \tau )  \in \cU_A \times \cC_A   \times \cM(A) $. To prove the continuity of $\Psi$, it suffices to show that $\phi$ and $\psi$ are continuous.  Hence we only need to show that $\phi$ is continuous.
For any $(u', B',\tau'), (\widetilde u, \widetilde B, \widetilde \tau )  \in \cU_A \times \cC_A   \times \cM(A) $, we have
\begin{align*}
\phi(u',B', \tau') - \phi (\widetilde u, \widetilde B, \widetilde \tau) &= \mathop{\max}\limits_{  a' \in  B' }  u' ( a',  \tau' )   -   \mathop{\max}\limits_{  \widetilde a \in \widetilde B } \widetilde u (\widetilde a, \widetilde \tau ) \\
&	\le \mathop{\max}\limits_{  a' \in  B' }   \bigl\{u' ( a',  \tau' )  -  u' ( a', \widetilde \tau )  \bigr\}  \tag{i$'$} \\
& + \mathop{\max}\limits_{  a' \in  B' }   \bigl\{u' ( a',  \widetilde \tau )  -  \widetilde u ( a', \widetilde \tau )  \bigr\}  \tag{ii$'$} \\
& + \mathop{\max}\limits_{  a' \in  B' } \widetilde u ( a', \widetilde \tau )  - \mathop{\max}\limits_{  \widetilde a \in \widetilde B } \widetilde u (\widetilde a, \widetilde \tau ).  \tag{iii$'$}
\end{align*}
Since  $u'$ is uniformly continuous on $A \times \cM(A)$, part (i$'$) tends to $0$ as $\rho(\tau', \widetilde \tau  ) \to 0$. Part (ii$'$) also tends to $0$ as $ ||u'- \widetilde u ||_{\infty} \to 0 $. Since $\widetilde u( \cdot, \widetilde \tau )$ is bounded and continuous on the compact set $A$, we can see that $ \Phi(\widetilde B) =  \max_{  \widetilde a \in \widetilde B } \widetilde u (\widetilde a, \widetilde \tau )
$ is also continuous function on $\cC_A$, where $\cC_A$ is endowed with the Hausdorff metric $d_H$. Hence part (iii$'$) tends to $0$ as $d_{H}( B', \widetilde B  ) \to 0$, which implies that $\Psi$ is a continuous function on $\cC_A \times \mathcal{U}_A \times \mathcal{M}(A) \times \mathcal{M}( A) $.

\bigskip

\noindent\textbf{Step 3.} Now we can finish the proof based on previous two steps. By simple calculations, we have that
\begin{align*}
&\min\limits_{a^{\prime} \in A_i^n }  \left\{ \int_A u_i^n \bigl(a,s(g^n) \bigr) g^n(i,\rmd a) - u_i^n \bigl(a^{\prime},s(a^{\prime},  g^n_{- i}  )  \bigr) \right\} \\
= &\min\limits_{a^{\prime} \in A_i^n }  \left\{ \int_A u_i^n \bigl(a,s(g^n)\bigr)  g^n(i,\rmd a) - u_i^n\bigl(a^{\prime},s(g^n)\bigr) + u_i^n\bigl(a^{\prime},s(g^n)\bigr) - u_i^n\bigl(a^{\prime}, s(a^{\prime},  g^n_{-i}) \bigr) \right\} \\
 \le &\min\limits_{a^{\prime} \in A_i^n }  \left\{ \int_A u_i^n \bigl(a,s(g^n)\bigr)  g^n(i,\rmd a) - u_i^n\bigl(a^{\prime},s(g^n)\bigr) \right\}
 + \max\limits_{a^{\prime} \in A_i^n }  \left\{  u_i^n\bigl(a^{\prime},s(g^n)\bigr) - u_i^n\bigl(a^{\prime}, s(a^{\prime},  g^n_{-i})\bigr)  \right\} \\
= &\Psi \bigl( G^n(i), g^n(i), s(g^n) \bigr) + \max\limits_{a^{\prime} \in A_i^n }  \left\{  u_i^n\bigl(a^{\prime},s(g^n)\bigr) - u_i^n\bigl(a^{\prime}, s(a^{\prime},  g^n_{-i}) \bigr) \right\}.
\end{align*}
By Lemma~\ref{lem-payoff gap}, we know that
$$\min_{a^{\prime} \in A_i^n }  \left\{ \int_A u_i^n \bigl(a,s(g^n)\bigr) g^n(i,\rmd a) - u_i^n \bigl(a^{\prime},s(a^{\prime},  g^n_{- i}  )\bigr)\right\}   \ge -\frac{\varepsilon}{2}$$
for all $i \in S^n$, $n \ge N_{\varepsilon}$.
Based on the Inequality~(\ref{equa-metric gap 2}) and equicontinuity of payoff functions, we know that there exists $\overline N \ge N_{\varepsilon}$ such that 
$$\max_{a^{\prime} \in A_i^n }  \left\{  u_i^n\bigl(a^{\prime},s(g^n)\bigr) - u_i^n\bigl(a^{\prime}, s(a^{\prime},  g^n_{-i})\bigr) \right\} \le \frac{\varepsilon}{2}$$ 
for all $i \in S^n, n \ge \overline N.$ Therefore, we have 
$$\Psi \bigl( G^n(i), g^n(i), s(g^n)\bigr) \geq -{\varepsilon}$$ 
for all $i^n \in S^n , n \ge \overline N.$

Let $h^n = \Psi\bigl(G^n,g^n,s(g^n)\bigr), \widetilde h = \Psi\bigl(G,g,s(g)\bigr)$. Since $\Psi$ is continuous and $\Bigl\{\cS \bigl(G^n,g^n,s(g^n)\bigr)\Bigr\}_{ n \in \mathbb{Z}_{+}  }$ converges weakly to $ \cS\bigl(G,g,s(g)\bigr)$, we conclude that $\{ h^n\}_{n \in \bZ_{+}}$ also converges weakly to $\widetilde h$. Thus we have
\[ 1 - \gamma \le \limsup\limits_{n \to \infty }  \lambda^n(h^n)^{ - 1}\bigl([ - \varepsilon ,\infty )\bigr) \le \lambda {\widetilde h}^{ - 1} \bigl([ -\varepsilon ,\infty )\bigr).\]
Letting $\varepsilon \to 0$ first and then letting $\gamma \to 0$, we
have ${\lambda {\widetilde h}^{-1}([0,\infty))}  = 1 $, which implies that 
$$\int_A u_i\bigl(a,s(g)\bigr) g( i, \rmd a) \geq u_i\bigl(a^{\prime},s(g)\bigr)$$
for $\lambda$-almost all $i \in I$ with $ a' \in A_i$.
Therefore,  $g$ is a Nash equilibrium of $G$.

%%%%%%%%%%%%%%%%%%%%%%%%%%%%%%%%%%%%%%%%%%%%%%%%%%%%%%%%%%%%%%%%%%%%%%%%%%%%%%%%%%%%%%%

\subsection{Proof of Proposition~\ref{prop-strong}}\label{subsec-proof-strong}
By Step 1 of the proof of Theorem~\ref{thm-closed graph property}, we know that  $\mathrm{supp\,\,} g(i)  \subseteq A_i$ for $\lambda$-almost all $i \in I$. Since $g^n $ is a Nash equilibrium of the large game $G^n$, we have
\[
u^n_i \bigl( g^n(i), s(g^n) \bigr)
\ge
u^n_i \bigl(a, s(g^n)\bigr) 
\]
for all $ i \in I^n$  and $a \in A_i^{n}$. 

By Step 2 of the proof of Theorem~\ref{thm-closed graph property}, we know that the   mapping $\Psi \colon \cC_A \times \cU_A \times \cM(A) \times \cM(A) \to \bR$:
 $$\Psi   (\widetilde{B}, \widetilde u, \widetilde \mu, \widetilde \tau) = \min_{a' \in \widetilde{B}} \left\{\int_A {\widetilde u(a, \widetilde\tau)} \widetilde \mu(\rmd a) - \widetilde u(a', \widetilde \tau)\right\}$$ 
 for any $(\widetilde{B}, \widetilde u, \widetilde \mu, \widetilde \tau )  \in  \cC_A \times \cU_A \times \cM(A) \times \cM(A) $,  is a continuous function on $\cC_A \times \cU_A \times \cM(A) \times \cM(A)$. By the definition of $\Psi$, we have that
\begin{align*}
	&\min\limits_{a^{\prime} \in A_i^n }  \left\{ \int_A u_i^n \bigl(a,s(g^n) \bigr) g^n(i,\rmd a) - u_i^n \bigl(a^{\prime},s( g^n  )  \bigr) \right\} \\
	= &\Psi \bigl( G^n(i), g^n(i), s(g^n) \bigr). 
\end{align*}
Since $g^n$ is a Nash equilibrium of game $G^n$, we know that
$$\min_{a^{\prime} \in A_i^n }  \left\{ \int_A u_i^n \bigl(a,s(g^n)\bigr) g^n(i,\rmd a) - u_i^n \bigl(a^{\prime},s(  g^n )\bigr)\right\}   \ge 0$$
for $\lambda$-almost all $i \in I^n$, for each $n \in \bN$.
Therefore, we have 
$$\Psi \bigl( G^n(i), g^n(i), s(g^n)\bigr) \geq 0$$ 
for $\lambda$-almost all $i \in I^n$, for each $n \in \bN$.

Next, we adopt $h$ and $\widetilde h$ as defined in Step 3 of the proof of Theorem~\ref{thm-closed graph property}.
Let $h^n = \Psi\bigl(G^n,g^n,s(g^n)\bigr), \widetilde h = \Psi\bigl(G,g,s(g)\bigr)$. Since $\Psi$ is continuous and $\Bigl\{\cS \bigl(G^n,g^n,s(g^n)\bigr)\Bigr\}_{ n \in \mathbb{Z}_{+}  }$ converges weakly to $ \cS\bigl(G,g,s(g)\bigr)$, we conclude that $\{ h^n\}_{n \in \bZ_{+}}$ also converges weakly to $\widetilde h$. Thus we have
\[ 1 \le \limsup\limits_{n \to \infty }  \lambda^n(h^n)^{ - 1}\bigl([0 ,\infty )\bigr) \le \lambda {\widetilde h}^{ - 1} \bigl([ 0 ,\infty )\bigr),\]
which implies that 
$$\int_A u_i\bigl(a,s(g)\bigr) g( i, \rmd a) \geq u_i\bigl(a^{\prime},s(g)\bigr)$$
for $\lambda$-almost all $i \in I$ with $ a' \in A_i$.
Therefore,  $g$ is a Nash equilibrium of $G$.

%%%%%%%%%%%%%%%%%%%%%%%%%%%%%%%%%%%%%%%%%%%%%%%%%%%%%%%%%%%%%%%%%%%%%%%%%%%%%%%%%%%%%%%

\subsection{Proof of Proposition~\ref{prop-refine}}\label{subsec-proof-pure}
This proof relies on the concept of the Nash equilibrium distribution of a large game. We begin by stating the definition of the Nash equilibrium distribution as follows.

\begin{defn}\label{defn-NED}
     \rm
	A probability measure $\zeta$ on $\cC_A \times \mathcal{U}_A \times A$ is a Nash equilibrium  distribution of a large game $G$ if\\
	(i) $\zeta|_{\cC_A \times \cU_A}=\lambda G^{-1};$ \\
	(ii) $\zeta ( \mathrm{Br}( \zeta) )=1$ where $\mathrm{Br}(\zeta )=\{(B, u ,a)\in \cC_A \times \cU_A \times A| u(a,\zeta|_A)\geq u (\widetilde a,\zeta|_{A})\,\,\textrm{for all}\,\, \widetilde a \in B \}$.
The set of Nash equilibrium distributions of a large game $G$ is denoted by $\mathrm{NED}(G)$.
\end{defn}

By Theorem~\ref{thm-closed graph property}, we obtain a symmetric Nash equilibrium $g$ that induces the distribution $\tau$, which is the limit of $\{ \int_{I^n} g^n \rmd \lambda^n \}_{ n \in \bZ_{+} }$. It is direct to check that $s(G,g) = \int_I  \delta_{G(i)}\otimes g(i) {\rm{d}}\lambda(i)  \in \text{NED}(G)$.
\begin{comment}
 By \citet[Theorem 3]{SSY2020}, there exists a symmetric Nash equilibrium $g^s$ such that
$s(G,g^s) = s(G,g)$ as well as $s(g^s) = s(g)$. Hence, the symmetric closed graph property holds.
\end{comment}
Let $\sigma(G)$ be the $\sigma$-algebra generated by the large game $G$.
If $\mathcal{I}$ is nowhere equivalent to $\sigma(G)$, then there exists a measurable mapping $f$ such that $s(G,g) =  \lambda (G,f)^{-1}$ (\citet[Lemma 2]{HSS2017}).
%Following
%the same argument in step 1 and step 2 of the proof of Theorem~\ref{thm-closed graph property} , we can obtain that $f$ is a  pure-strategy Nash equilibrium of $G$.
We divide the remaining proof into two steps. In step 1, we show that
$f(i) \in A_i$ for $\lambda$-almost all $i \in I$. In step 2, we show that $f$ is a Nash equilibrium.

\noindent\textbf{Step 1.} Since  $\{ s( \cA^n , g^n  ) \}_{n \in \bZ_{+}}$ converges weakly to $s( \cA , g)$, where  $s( \cA^n , g^n  ) = \int_{I^n}  \delta_{\cA^n(j)} \otimes g^n(j) \rmd \lambda^n(j)$ and $s(\cA,g) = \int_I  \delta_{\cA(j)} \otimes g(j) \rmd \lambda(j) $,
  $\{ s( \cA^n , g^n  ) \}_{n \in \bZ_{+}}$ also converges weakly to $\lambda ( \cA , f)^{-1}$.
Let $Z = \{ (B,b) | B \in \cC_A, b \in B \}$. Then $Z$ is a closed subset and
$s(\cA^n, g^n)(Z) = 1 $ for all $n \in \mathbb{Z}_{+}$. By weak convergence, we have $\lambda (\cA, f)^{-1}(Z) = 1.$ Since  $\lambda (\cA,f)^{-1}|_{\cC_A}(\cA(I)) = 1$, we conclude that $f(i) \in A_i$ for $\lambda$-almost all $i \in I$.

\noindent\textbf{Step 2.}  Since $\nu = s(G,g)= \lambda(G,f)^{-1} \in \mathrm{NED}(G)$, we have $\nu ( \text{Br}(\nu) ) = \lambda ( \{ i \in I  \colon (G(i), f(i) ) \in \text{Br}(\nu)     \})  = 1$. By the definition of $\text{Br}(\nu)$ and the fact that $\nu|_{ A} = \lambda f^{-1}$, we know that  for $\lambda$-almost all $i \in I$, we have $u_i(f(i),\lambda f^{-1})\geq u_i (\widetilde a, \lambda f^{-1})$  for all $\widetilde a \in A_i $. Thus, $f$ is a pure strategy Nash equilibrium of $G$ and $\lambda f^{-1} = s(g) = \tau$.

%------------------------------------------------------------------

{\small
\singlespacing

\end{document}